# Numerical Solution of a Singularly –Perturbed Boundary-Value Problems by Using A Non-Polynomial Spline

## ISLAM KHAN[1] AND TARIQ AZIZ[2]

[1]Department of Mathematics, Yanbu University College,

Yanbu Al-Sinaiyah, P.O. Box 31387, Kingdom of Saudi Arabia

[2]Department of Applied Mathematics, Faculty of Engg. & Technology,

Aligarh Muslim University, Aligarh−202002 (U.P.), India

**Abstract:**

We consider a non-polynomial cubic spline to develop the classes of methods for the numerical solution of singularly perturbed two-point boundary value problems. The proposed methods are second and fourth order accurate and applicable to problems both in singular and non-singular cases. Numerical results are given to illustrate the efficiency of our methods and compared with the methods given by different authors.

## 1. INTRODUCTION

We consider a second−order singularly perturbed boundary value problem

$$-\varepsilon y'' + P(x)y = f(x) \tag{1}$$

$$y(0)=A, \quad y(1)=B, \tag{2}$$

---

*Corresponding author.

E−mail: isk003@yahoo.co.in (I. Khan)

where A and B are given constants and ε is a parameter such that 0<ε<<1 and P(x), f(x) are smooth, bounded, real functions P:R→R, and f:R→R. The approximate solution of boundary−value problems with a small parameter affecting highest derivative of the differential equation is described. It is a well−known fact that the solution of singularly perturbed boundary−value problem exhibits a multiscale character. That is, there is a thin layer where the solution varies rapidly, while away from the layer the solution behaves regularly and varies slowly. This class of problems has recently gained importance in the literature for two main reasons. Firstly, they occur frequently in many areas of science and engineering, for example, combustion, chemical reactor theory, nuclear engineering, control theory, elasticity, fluid mechanics etc. A few notable examples are boundary−layer problems, WKB Theory, the modelling of steady and unsteady viscous flow problems with large Reynolds number and convective heat transport problems with large Peclet number. Secondly, the occurrence of sharp boundary−layers as ε, the coefficient of highest derivative, approaches zero creates difficulty for most standard numerical schemes. The application of spline for the numerical solution of singularly-perturbed boundary value-problems has been described by many authors [2-5, 7-11].

In the present paper, we have derived a three point formula based on non-polynomial cubic spline. Analysis of the method shows that it has second−order convergence for arbitrary $\lambda_1, \lambda_2$ such that $\lambda_1+\lambda_2=\frac{1}{2}$ and fourth−order convergence for $\lambda_1=\frac{1}{12}, \lambda_2=\frac{5}{12}$. In section 2, we have given derivation of the method and in section 3, Truncation error. Section 4 contains numerical illustrations and application of methods on two examples is given.

## 2. DERIVATION OF THE METHOD

Let $x_0=a,\ x_N=b,\ x_i=a+ih,\ h=(b-a)/N$.

A function S($x,\tau$) of class C$^2$[a,b] which interpolates y($x$) at the mesh point $x_i$ depends on a parameter $\tau$, reduces to cubic spline in [a,b] as $\tau \to 0$ is termed as parametric cubic spline function.

The spline function S($x,\tau$) = S($x$) satisfying in [$x_i, x_{i+1}$], the differential equation

$$S''(x) - \tau S(x) = [S''(x_i) - \tau S(x_i)]\frac{(x_{i+1} - x)}{h} + [S''(x_{i+1}) - \tau S(x_{i+1})]\frac{(x - x_i)}{h}, \qquad (3)$$

where S($x_i$) = $y_i$ and $\tau > 0$ is termed as cubic spline in tension.

Solving the linear second−order differential equation (3) and determining the arbitrary constants from the interpolatory conditions S($x_{i+1}$) = $y_{i+1}$, S($x_i$)=$y_i$ we get, after writing $\lambda = h\tau^{1/2}$

$$S(x) = \frac{h^2}{\lambda^2 \sinh \lambda}\left[ M_{i+1} \sinh \frac{\lambda(x - x_i)}{h} + M_i \sinh \frac{\lambda(x_{i+1} - x)}{h} \right]$$

$$- \frac{h^2}{\lambda^2}\left[ \frac{(x - x_i)}{h}\left( M_{i+1} - \frac{\lambda^2}{h^2} y_{i+1} \right) + \frac{(x_{i+1} - x)}{h}\left( M_i - \frac{\lambda^2}{h^2} y_i \right) \right]. \qquad (4)$$

Differentiating equation (4) and letting $x \to x_i$, we obtain

$$S'(x_i^+) = \frac{y_{i+1} - y_i}{h} - \frac{h}{\lambda^2}\left[ \left( 1 - \frac{\lambda}{\sinh \lambda} \right) M_{i+1} + (\lambda \coth \lambda - 1) M_i \right].$$

Considering the interval ($x_{i-1}, x_i$) and proceeding similarly, we obtain

$$S'(x_i^-) = \frac{y_i - y_{i-1}}{h} + \frac{h}{\lambda^2}\left[ (\lambda \coth \lambda - 1) M_i + \left( 1 - \frac{\lambda}{\sinh \lambda} \right) M_{i-1} \right].$$

Equating the left− and right−hand derivatives at $x_i$, we have

$$\frac{y_i - y_{i-1}}{h} + \frac{h}{\lambda^2}\left[ (\lambda \coth \lambda - 1) M_i + \left( 1 - \frac{\lambda}{\sinh \lambda} \right) M_{i-1} \right]$$

$$= \frac{y_{i+1} - y_i}{h} - \frac{h}{\lambda^2}\left[\left(1 - \frac{\lambda}{\sinh \lambda}\right)M_{i+1} + (\lambda \coth \lambda - 1)M_i\right]. \qquad (5)$$

This leads to the tridiagonal system

$$h^2(\lambda_1 M_{i-1} + 2\lambda_2 M_i + \lambda_1 M_{i+1}) = y_{i+1} - 2y_i + y_{i-1}, \qquad (6)$$

where $\lambda_1 = \frac{1}{\lambda^2}\left(1 - \frac{\lambda}{\sinh \lambda}\right)$, $\lambda_2 = \frac{1}{\lambda^2}(\lambda \coth \lambda - 1)$, $M_i = S''(x_i), i = 1(1)N - 1$.

For a numerical solution of the boundary value problem (1) the interval [0, 1] is divided into a set of grid points with step length $h = (b-a)/N$, N being a positive integer the spline approximation on [0, 1] that consist of the nodal points $x_0 = a$, $x_N = b$, $x_i = a + ih$, $i = 1(1)N - 1$, then

$$-\varepsilon y''(x_i) + P(x_i) y(x_i) = f(x_i) \qquad (7)$$

or $\quad y''(x_i) = \frac{P(x_i)}{\varepsilon} y(x_i) - \frac{1}{\varepsilon} f(x_i)$

by using the spline's second derivative we have

$$\begin{cases} M_{i-1} = \dfrac{P_{i-1}}{\varepsilon} y_{i-1} - \dfrac{1}{\varepsilon} f_{i-1} \\[2mm] M_i = \dfrac{P_i}{\varepsilon} y_i - \dfrac{1}{\varepsilon} f_i \\[2mm] M_{i+1} = \dfrac{P_{i+1}}{\varepsilon} y_{i+1} - \dfrac{1}{\varepsilon} f_{i+1} \end{cases} \qquad (8)$$

Substitute (8) in (6) we have

$$(\lambda_1 h^2 P_{i-1} - \varepsilon) y_{i-1} + 2(\varepsilon + \lambda_2 h^2 P_i) y_i + (\lambda_1 h^2 P_{i+1} - \varepsilon) y_{i+1}$$
$$= h^2(\lambda_1 f_{i-1} + 2\lambda_2 f_i + \lambda_1 f_{i+1}), \qquad i = 1(1)N - 1. \qquad (9)$$

Finally we have the following system

$$2(\varepsilon + \lambda_2 h^2 P_1) y_1 + (\lambda_1 h^2 P_2 - \varepsilon) y_2$$
$$= h^2(\lambda_1 f_0 + 2\lambda_2 f_1 + \lambda_1 f_2) y_{i-1} - (\lambda_1 h^2 P_0 - \varepsilon)A, \qquad i = 1$$

$$(\lambda_1 h^2 P_{i-1} - \varepsilon) y_{i-1} + 2(\varepsilon + \lambda_2 h^2 P_i) y_i + (\lambda_1 h^2 P_{i+1} - \varepsilon) y_{i+1}$$
$$= h^2(\lambda_1 f_{i-1} + 2\lambda_2 f_i + \lambda_1 f_{i+1}), \qquad 2 \leq i \leq N-2 \qquad (10)$$
$$(\lambda_1 h^2 P_{N-2} - \varepsilon) y_{N-2} + 2(\varepsilon + \lambda_2 h^2 P_{N-1}) y_{N-1}$$
$$= h^2(\lambda_1 f_{N-2} + 2\lambda_2 f_{N-1} + \lambda_1 f_N) - (\lambda_1 h^2 P_N - \varepsilon) B, \qquad i = N-1 \qquad (11)$$

## 3. Truncation error

By expanding equation (9) in Taylor's series about $x_i$ we obtain

$$T_i(h) = [-1 + 2(\lambda_1 + \lambda_2)] \varepsilon h^2 y^{(2)}(x_i) + (\lambda_1 - \frac{1}{12}) \varepsilon h^4 y^{(4)}(x_i) + [\frac{\lambda_1}{12} - \frac{1}{360}] \varepsilon h^6 y^{(6)}(x_i) + O(h^8)$$

by choosing different values of $\lambda_1$ and $\lambda_2$ provided $\lambda_1 + \lambda_2 = \frac{1}{2}$, we can obtain the classes of second order methods.

**Remark 1:** if $\lambda_1 = \frac{1}{6}$ and $\lambda_2 = \frac{1}{3}$, $T_i = O(h^4)$ then the resulting method is reduced to the cubic spline that is second order method.

**Remark 2:** if $\lambda_1 = \frac{1}{12}$ and $\lambda_2 = \frac{5}{12}$, $T_i = O(h^6)$ then the resulting method is fourth- order method.

## 4. Numerical illustrations and discussion

We have implemented our method on two examples of singularly perturbed boundary value problems. These examples have been solved for $\lambda_1 = \frac{1}{12}$ and $\lambda_2 = \frac{5}{12}$ and different values of h. The numerical solutions are computed and compared with the exact solutions at grade points. The maximum absolute errors $E = \max \left| \bar{y}_i - y_i \right|$ in numerical solutions are tabulated in table 1 and 2.

**Example 4.1**

$$-\varepsilon\, y'' + y = -\cos^2(\pi x) - 2\varepsilon\pi^2\cos(2\pi x),$$

$$y(0) = y(1) = 0.$$

The exact solution is given by

$$y(x) = [\exp(-(1-x)/\sqrt{\varepsilon}) + \exp(-x/\sqrt{\varepsilon})]/[1 + \exp(-1/\sqrt{\varepsilon})] - \cos^2(\pi x),$$

This problem has been solved using the method (9) with different values of $N = 16,\ 32,\ 64,\ 128,\ 256$ and $\varepsilon = \dfrac{1}{16}, \dfrac{1}{132}, \dots \dfrac{1}{128}$. The maximum absolute errors in solutions are tabulated in table 1 and compared with the results in [8-11], which shows the accuracy of our methods.

**Example 4.2**

$$-\varepsilon\, y'' + (1+x)\, y = -40[x(x^2 - 1) - 2\varepsilon]$$

$$y(0) = y(1) = 0.$$

The exact solution is given by $y(x) = 40\, x(1-x)$.

This problem has been solved using the method (9) with different values of $N = 16,\ 32$ and $\varepsilon = 0.1 \times 10^{-3}, \dots \dots 0.1 \times 10^{-8}$. The maximum absolute errors in solutions are tabulated in table 2 and compared with the results in [3, 5, 7] which shows the superiority of our methods.

**Table 1: Maximum absolute errors in solutions of problem 1**

**Our method**

| $\varepsilon$ | N=16 | N=32 | N=64 | N=128 | N=256 |
|---|---|---|---|---|---|
| 1/16 | $6.09 \times 10^{-7}$ | $0.780 \times 10^{-8}$ | $1.32 \times 10^{-7}$ | $9.98 \times 0^{-9}$ | $1.19 \times 10^{-15}$ |
| 1/32 | $1.12 \times 10^{-6}$ | $1.24 \times 10^{-8}$ | $8.87 \times 10^{-8}$ | $6.52 \times 10^{-9}$ | $4.62 \times 10^{-15}$ |
| 1/64 | $3.54 \times 10^{-6}$ | $2.78 \times 10^{-7}$ | $7.89 \times 10^{-7}$ | $2.54 \times 10^{-8}$ | $9.14 \times 10^{-10}$ |
| 1/128 | $2.27 \times 10^{-5}$ | $1.23 \times 10^{-7}$ | $5.41 \times 10^{-7}$ | $5.55 \times 10^{-8}$ | $4.78 \times 10^{-9}$ |

**Surja and Staojanovic method [9]**

| $\varepsilon$ | N=16 | N=32 | N=64 | N=128 | N=256 |
|---|---|---|---|---|---|
| 1/16 | $8.06 \times 10^{-3}$ | $2.02 \times 10^{-3}$ | $5.08 \times 10^{-4}$ | $1.27 \times 10^{-4}$ | $3.17 \times 10^{-5}$ |
| 1/32 | $7.11 \times 10^{-3}$ | $1.79 \times 10^{-3}$ | $4.48 \times 10^{-4}$ | $1.12 \times 10^{-4}$ | $2.80 \times 10^{-5}$ |
| 1/64 | $6.58 \times 10^{-3}$ | $1.66 \times 10^{-3}$ | $4.15 \times 10^{-4}$ | $1.04 \times 10^{-4}$ | $2.60 \times 10^{-5}$ |
| 1/128 | $6.36 \times 10^{-3}$ | $1.61 \times 10^{-3}$ | $4.03 \times 10^{-4}$ | $1.01 \times 10^{-4}$ | $2.52 \times 10^{-5}$ |

**Surja and Herceg and Cvetkovic's method [10]**

| $\varepsilon$ | N=16 | N=32 | N=64 | N=128 | N=256 |
|---|---|---|---|---|---|
| 1/16 | $4.14 \times 10^{-4}$ | $1.02 \times 10^{-3}$ | $2.54 \times 10^{-4}$ | $6.35 \times 10^{-5}$ | $1.58 \times 10^{-5}$ |
| 1/32 | $3.68 \times 10^{-3}$ | $9.03 \times 10^{-4}$ | $5.61 \times 10^{-5}$ | $1.40 \times 10^{-5}$ | $3.50 \times 10^{-5}$ |
| 1/64 | $3.45 \times 10^{-3}$ | $8.40 \times 10^{-4}$ | $2.08 \times 10^{-4}$ | $5.20 \times 10^{-5}$ | $1.30 \times 10^{-5}$ |
| 1/128 | $3.43 \times 10^{-3}$ | $8.21 \times 10^{-4}$ | $2.03 \times 10^{-4}$ | $5.06 \times 10^{-5}$ | $1.26 \times 10^{-5}$ |

**Surja and  Vukoslavcevic's method [11]**

| $\varepsilon$ | N=16 | N=32 | N=64 | N=128 | N=256 |
|---|---|---|---|---|---|
| 1/16 | $1.20 \times 10^{-4}$ | $7.47 \times 10^{-6}$ | $4.67 \times 10^{-7}$ | $2.90 \times 10^{-8}$ | $4.39 \times 10^{-9}$ |
| 1/32 | $1.28 \times 10^{-4}$ | $8.0 \times 10^{-6}$ | $5.0 \times 10^{-7}$ | $3.14 \times 10^{-8}$ | $1.99 \times 10^{-9}$ |
| 1/64 | $1.60 \times 10^{-4}$ | $1.0 \times 10^{-5}$ | $6.26 \times 10^{-7}$ | $3.92 \times 10^{-8}$ | $2.31 \times 10^{-9}$ |
| 1/128 | $2.34 \times 10^{-4}$ | $1.47 \times 10^{-5}$ | $9.23 \times 10^{-7}$ | $5.77 \times 10^{-8}$ | $3.72 \times 10^{-9}$ |

**Kadalbajoo and  Bawa's method [8]**

| $\varepsilon$ | N=16 | N=32 | N=64 | N=128 | N=256 |
|---|---|---|---|---|---|
| 1/16 | $7.09 \times 10^{-3}$ | $1.77 \times 10^{-3}$ | $4.45 \times 10^{-4}$ | $1.11 \times 10^{-4}$ | $2.78 \times 10^{-5}$ |
| 1/32 | $5.68 \times 10^{-3}$ | $1.42 \times 10^{-3}$ | $3.55 \times 10^{-4}$ | $8.89 \times 10^{-5}$ | $2.22 \times 10^{-5}$ |
| 1/64 | $4.07 \times 10^{-3}$ | $1.01 \times 10^{-3}$ | $2.54 \times 10^{-4}$ | $6.35 \times 10^{-5}$ | $1.58 \times 10^{-5}$ |
| 1/128 | $6.97 \times 10^{-3}$ | $1.75 \times 10^{-3}$ | $4.33 \times 10^{-4}$ | $1.08 \times 10^{-4}$ | $2.71 \times 10^{-5}$ |

**Table 2: Maximum absolute errors in solutions of problem 2**

**For N=16**

| $\varepsilon$ | Method in [3] | Method in [5] | Method in [7] | Our method |
|---|---|---|---|---|
| $0.1 \times 10^{-3}$ | $0.25 \times 10^{-1}$ | $0.26 \times 10^{-1}$ | $0.65 \times 10^{-4}$ | $0.78 \times 10^{-15}$ |
| $0.1 \times 10^{-4}$ | $0.21 \times 10^{-1}$ | $0.24 \times 10^{-1}$ | $0.36 \times 10^{-4}$ | $0.76 \times 10^{-15}$ |
| $0.1 \times 10^{-5}$ | $0.70 \times 10^{-2}$ | $0.17 \times 10^{-1}$ | $0.33 \times 10^{-4}$ | $0.87 \times 10^{-15}$ |
| $0.1 \times 10^{-6}$ | $0.75 \times 10^{-3}$ | $0.69 \times 10^{-2}$ | $0.26 \times 10^{-4}$ | $0.91 \times 10^{-15}$ |
| $0.1 \times 10^{-7}$ | $0.74 \times 10^{-4}$ | $0.23 \times 10^{-2}$ | $0.20 \times 10^{-4}$ | $0.65 \times 10^{-15}$ |
| $0.1 \times 10^{-8}$ | $0.67 \times 10^{-5}$ | $0.76 \times 10^{-3}$ | $0.11 \times 10^{-4}$ | $2.71 \times 10^{-15}$ |

**For N=32**

| $\varepsilon$ | Method in [3] | Method in [5] | Method in [7] | Our method |
|---|---|---|---|---|
| $0.1 \times 10^{-3}$ | $0.64 \times 10^{-2}$ | $0.65 \times 10^{-2}$ | $0.59 \times 10^{-4}$ | $1.28 \times 10^{-15}$ |
| $0.1 \times 10^{-4}$ | $0.61 \times 10^{-2}$ | $0.64 \times 10^{-2}$ | $0.21 \times 10^{-4}$ | $1.36 \times 10^{-15}$ |
| $0.1 \times 10^{-5}$ | $0.41 \times 10^{-2}$ | $0.56 \times 10^{-2}$ | $0.35 \times 10^{-4}$ | $1.36 \times 10^{-15}$ |
| $0.1 \times 10^{-6}$ | $0.77 \times 10^{-3}$ | $0.31 \times 10^{-2}$ | $0.39 \times 10^{-4}$ | $1.49 \times 10^{-15}$ |
| $0.1 \times 10^{-7}$ | $0.76 \times 10^{-4}$ | $0.12 \times 10^{-2}$ | $0.21 \times 10^{-4}$ | $1.65 \times 10^{-15}$ |
| $0.1 \times 10^{-8}$ | $0.67 \times 10^{-5}$ | $0.38 \times 10^{-3}$ | $0.21 \times 10^{-4}$ | $3.71 \times 10^{-15}$ |